\documentclass[12pt,reqno]{amsart}
\usepackage{amsmath,amssymb,latexsym,textcomp,mathrsfs}
\usepackage[all]{xy}
\usepackage{graphicx}
\usepackage{bm,amsmath, amsthm, amssymb, amsfonts}
\usepackage{times}
\usepackage[english]{babel}

\setbox0=\hbox{$+$}
\newdimen\plusheight
\plusheight=\ht0
\def\+{\;\lower\plusheight\hbox{$+$}\;}

\setbox0=\hbox{$-$}
\newdimen\minusheight
\minusheight=\ht0
\def\-{\;\lower\minusheight\hbox{$-$}\;}

\setbox0=\hbox{$\cdots$}
\newdimen\cdotsheight
\cdotsheight=\plusheight%\ht0
\def\cds{\lower\cdotsheight\hbox{$\cdots$}}

\numberwithin{equation}{section}
\theoremstyle{plain}

  \newenvironment{nouppercase}{%
   \renewcommand{\uppercasenonmath}[1]{}}{}
	
	 \newcommand{\Keywords}[1]{\par\noindent
   {\small{Keywords and phrases}: #1}}
   
   \newcommand{\AMS}[1]{\par\noindent
   {\small{AMS Subject Classification (2010)}: #1}}

\begin{document}

\title{$\lambda^*$-CLOSED SETS AND  NEW SEPARATION AXIOMS IN ALEXANDROFF SPACES}
 % \author{Amar Kumar Banerjee$^1$}
 %\author{Jagannath Pal$^2$}
 %\address{1. Department of Mathematics,The University of Burdwan,Golapbag,Burdwan-713104,West Bengal,India}
 %\email{akbanerjee1971@gmail.com}
 %\address{2. Department of Mathematics,The University of Burdwan,Golapbag,Burdwan-713104,West Bengal,India}
 %\email{jpalbu1950@gmail.com}

 \author{Amar Kumar Banerjee$^1$}
 \author{Jagannath Pal$^2$}
 \newcommand{\acr}{\newline\indent}
 \maketitle
 \address{{1\,} Department of Mathematics, The University of Burdwan, Golapbag, Burdwan-713104,
 West Bengal, India.
                 Email: akbanerjee@math.buruniv.ac.in,akbanerjee1971@gmail.com\acr
            {2\,} Department of Mathematics, The University of Burdwan, Golapbag, Burdwan-713104,
            West Bengal, India. 
            Email:jpalbu1950@gmail.com\\}   
\begin{abstract}
Here we have studied the ideas of $g^*$-closed sets,  $g\wedge_\tau $-sets and $\lambda^*$-closed sets and investigate some of their properties in the spaces of A. D. Alexandroff \cite{AD}. We have also studied few separation axioms like $ T_\frac{\omega}{4},  T_\frac{3\omega}{8},  T_\omega $ in Alexandroff spaces and have introduced a new separation axiom namely $ T_\frac{5\omega}{8} $ axiom in this space. 
\end{abstract}

\begin{nouppercase}
\maketitle
\end{nouppercase}

\let\thefootnote\relax\footnotetext{        
\AMS{ 54A05,   54D10}.
}

\section{\bf Introduction}
\label{sec:int}
Topological spaces have been generalised in many ways. A. D. Alexandroff (1940) weakened the union requirements where  only countable union of open sets were taken to be open.  

The idea of generalised closed sets in a topological space was given by Levine \cite{NL}. Later many works on generalised closed sets have been done \cite{AC}, \cite{CA}, \cite{DP}, \cite{WD} etc. In 2003, P.  Das et al obtained a generalisation of closed sets in Alexandroff spaces which was called $ g^* $-closed sets. They investigated various properties on $ g^*$-closed  sets and also introduced a new separation axiom namely $ T_w $ axiom in Alexandroff spaces exactly in the same way as Levine \cite{NL} defined $ T_\frac{1}{2} $-spaces in topological spaces. Recently M. S. Sarsak \cite{MS} studied the same in a generalised topological space where a generalised topology $ \mu $ on a nonempty set $ X $ is a collection of subsets of $ X $ such that $ \emptyset \in \mu $ and $ \mu $ is closed under arbitrary unions. Members of $ \mu $ are called $ \mu $-open sets. He also introduced some new separation axioms namely $ \mu$-$T_\frac{1}{4},  \mu$-$T_\frac{3}{8} $ and $ \mu$-$T_\frac{1}{2} $ axioms and studied their properties and relations between the axioms.

Here we have studied the idea of generalised closed sets namely the  $ g^*$-closed sets in Alexandroff spaces. We have investigated few more properties of $ g^*$-closed  sets which were not studied in \cite{DP}. Also we have obtained few separation axioms like $ T_\frac{\omega}{4},  T_\frac{3\omega}{8},  T_\omega $ in Alexandroff spaces in the same way as that of Sarsak \cite{MS} and investigate how far several results as valid in \cite{MS} are affected in Alexandroff spaces. Also we have introduced a new separation axiom namely $ T_\frac{5\omega}{8} $ axiom in a Alexandroff space which was not studied before. We have found out that $ T_\frac{5\omega}{8} $-space can be placed in between $ T_\frac{3\omega}{8} $ and $ T_\omega $ spaces.

 \section{\bf Preliminaries}
 \label{sec:pre}

\textbf{Definition 2.1}  \cite{AD}: A non empty set $X$ is called a $ \sigma $-space or simply a space if in it is chosen a system of subsets $\mathcal{F}$ satisfying the following axioms:

(1)	 The intersection of a countable number of sets in $\mathcal{F}$ is a set in $\mathcal{F}$.     

(2)  The union of a finite number of sets in $\mathcal{F}$ is a set in $\mathcal{F}$.

(3) The void set is a set in $\mathcal{F}$.

(4)	The whole set $ X $ is a set in $\mathcal{F}$.

Sets of $\mathcal{F}$ are called closed sets. Their complementary sets are called open sets. The collection of all such open sets will sometimes be denoted by $\tau $ and the space by    $(X ,\tau )$. When there is no confusion, the space  $(X , \tau)$ will simply be denoted by $ X $.

Note that a topological space is a space but in general $ \tau $ is not a topology as can be easily seen by taking $ X = R $ and $ \tau $ as the collection of all $F_\sigma$-sets in $ R $. Several examples of spaces are seen in \cite{PD}, \cite{DP}, \cite{LD}. The definition of closure of a set and interior of a set in a space are similar as in the case of a topological space. Note that closure of a set in a space may not be closed in general. Also interior of a set in a space may not be open.

Throughout the paper $ X $ stands for a space and   sets are always subsets of $ X $ unless otherwise stated. The letters $ R $ and $ Q $ stand respectively for the set of real numbers and the set of rational numbers.\\

\textbf{Definition 2.2} \cite{WD} :  Two sets $ A, B $ in $ X $ are said to be weakly separated if there are two open sets $U,  V $ such that $A\subset U, B\subset V $ and $ A\cap V=B\cap U=\emptyset$.\\

\textbf{Definition 2.3} \cite{AD} : A  space $(X, \tau)$ is called $T_0$ if for any pair of distinct points of $X$, there is an open set containing one of the points but not the other.

Observe that a space $(X,  \tau)$  is $T_0$ if and only if for any pair of distinct points $x, y \in X $,  there is a set $A$ containing one of the points such that $ A $ is either open or closed. 

Clearly it can be easily checked that if a space $ (X,  \tau )$ is $T_0$  then for every pair of distinct points $ p,q\in X $,  either $ p\not \in\overline{\{q\}} $  or $ q\not \in \overline
{\{p\}} $.\\

\textbf{Definition 2.4} \cite{AD} : $(X, \tau)$ is said to be $T_1$-space if for any pair of distinct points $x, y \in X$, there exist open sets  $ U, V $ such that $x \in U,  y \not \in U,  y\in V , x\not\in V$.\\

\textbf{Definition 2.5 } \cite{DP} : A point $ x\in X $ is said to be a limit point of  $ A $ in a space $ (X, \tau) $ if for any open set $ U $ containing $ x $,  $ U \cap(A - \{x\})\not=\emptyset $ . The set of all limit points of $ A $ is called derived set of $ A $ and is denoted by $ A' $. \\

\textbf{Definition  2.6} \cite{SM}  :  A space $ (X,  \tau) $ is said to be a $ R_0$-space  if every open set contains the closure of each of its singleton.\\

\textbf{Definition  2.7} \cite{AD} :  A space $ (X,  \tau) $ is said to be bicompact  if every open cover of it has a finite subcover.\\

\textbf{  Theorem  2.8 }:  Let $A$ be a subset in a space  $(X, \tau)$, then $ X-\overline{X-A}=Int(A) $.

The proof is straightforward and so is omitted.

\section{\bf    $\textit  g^*$-closed sets, \quad $\textit g\wedge_\tau $-sets \quad and \quad $\textit g\vee _\tau $-sets in a space}

\textbf{Definition 3.1 } \cite{NL}: A set $A $ in a topological space is said to be a generalised closed  ($ g $-closed for short) if and only if $\overline {A} \subset U$ whenever $A\subset U$ and $ U $ is open.\\

                             \textbf{Definition 3.2 :} (cf.  Definition 9 \cite {DP} ):  A set $ A $ in a space is said to be a $g^*$-closed set if and only if there is a closed set  $F$ containing $ A $ such that $ F \subset U$ whenever $A\subset U$ and  $U$ is open.
                             
                             Clearly every closed set is $g^*$-closed but the converse may not be true as shown in the Example 1 \cite{DP}.\\
                             
                             We now state few important theorems given in \cite{DP}.\\
                             
                                                         \textbf{Theorem  3.3} \cite{DP}:
 A set $A$ is $ g^*$-closed  if and only if there is a closed set $F$ containing $A$ such that $ F \subset$ ker$(A)$.\\

                             \textbf{Theorem 3.4} \cite{DP} :  For each $x \in X,$  $\{x\}$ is either closed or $ X - \{x\}$ is $ g^*$-closed.\\

                             \textbf{Theorem 3.5} \cite{DP} :  A set $A$ is $ g^*$-closed  if and only if there is a closed set $F$ containing $A$ such that $ F - A $ does not contain any non-void closed set.\\
                             
                             \textbf{Definition 3.6}  \cite{DP}: A set $A$ is called $ g^*$-open  set if $ X - A $ is $ g^*$-closed  set.\\

\textbf{Theorem 3.7} \cite{DP} :  A set $A$ is $ g^*$-open  set if and only if there is an open set $V$ contained in $A$ such that $ F \subset V $ whenever $F$ is closed and $ F \subset A $.\\

\textbf{Definition 3.8} (cf.\cite {DP}) :   Let $A \subset X$. The kernel of $ A$ is the set $\cap \{U: U \in \tau , A\subset U\}$ and denoted by $A_\tau^\wedge$.\\

\textbf{Definition 3.9 (cf. \cite {MS})} : Let $ A \subset X $. Then we define $A_\tau^\vee =  \cup \{F:  X - F \in \tau ,  F \subset A\}$.\\

\textbf{Definition 3.10} (cf. \cite {MS}):  A set $A$ is called a $\wedge_\tau$-set if $A = A_\tau^\wedge$ or, equivalently,  $A$ is the intersection of all open sets containing $A$. \\
 
 \textbf{Definition 3.11}  (cf. \cite{MS}):  A set $A$ is called a $\vee_\tau$-set if $A = A_\tau^\vee$ or, equivalently, $X - A$ is a $\wedge_\tau $-set (i.e. $A$ is the union of all closed sets contained in $A$). Clearly $ X $ and $ \emptyset $ are $ \vee_\tau $-sets as well as $\wedge_\tau$-sets.\\

Clearly an open set is $ g^*$-open  set. But the converse may not be true as shown in the following example.\\
  
\textbf{Example  3.12}: Let $ X = R - Q $,  $ \tau $ = $\{X,  \emptyset,  G_i\}$, Where $ G_i $ be the all countable subsets in $X$.  Thus $ (X, \tau) $ is a space but not a topological space.  Let $B$ be the set of all positive irrational numbers.  Then $B$ is $ g^*$-closed,  since $X$ is the only open and closed set containing $B$. Then $ X - B $ is $ g^*$-open   but is not an open set.\\

\textbf{Theorem 3.13 :  } Intersection of two $ g^* $-open sets is $ g^* $-open.

\begin{proof}
Let $A , B$ be two $ g^* $-open sets. There exist open sets $U , V$ contained in $A , B$ respectively such that $ F_1\subset U $ and $ F_2\subset V $ whenever $ F_1\subset A $ and $ F_2\subset B $;   $ F_1 $ and $ F_2 $ are closed sets. Clearly $ U\cap V $ is open and $ U\cap V \subset A\cap B $. Let $ F $ be a closed set such that $ F \subset A\cap B $. Then $ F\subset A $ and $ F\subset B $. So  $ F\subset U $ and $ F\subset V $. Therefore $ F\subset U\cap V $ which implies $ A\cap B $ is $ g^* $-open.
\end{proof}

\textbf{Theorem 3.14} :  A set $A$ is $ g^*$-open  if and only if there is an open set $V \subset A $ such that $ A_\tau^\vee \subset V $. 
\begin{proof}
Let $A$ be a $ g^*$-open set. Then by Theorem 3.7, there exists an open set $ V \subset A $ such that $ F \subset V $ whenever $ F \subset A $ and $ F $ is closed. So $ \cup \{F : F \subset A, F $ is closed\} $ \subset V $.  This implies that  $ A_\tau^\vee \subset V $.

Conversely, let there be an open set $ V \subset A $ such that $ A_\tau^\vee \subset V $. This implies that $ A_\tau^\vee = \cup \{F: F\subset A, F $ is closed\} $ \subset  V $. So there is an open set $ V $ such that $ F \subset V $ whenever $ F \subset A , F $ is closed. Hence $A$ is $ g^*$-open by Theorem 3.7.
\end{proof}

\textbf{Definition 3.15} (cf. Definition 2.11\cite{MS}) :  A set $A$ is called a generalised $\wedge_\tau$-set denoted as $g\wedge_\tau$-set if $A_\tau^\wedge\subset F$ whenever $F\supset A$ and $F$ is closed. A set $A$ is called a generalised $ \vee_\tau$-set  denoted as $ g\vee_\tau $-set if $ X - A  $ is $ g\wedge_\tau$-set.\\

\textbf{Theorem 3.16:} A set $A$ is a $ g\vee_\tau$-set  if and only if $ A_\tau^\vee\supset U $ whenever $A \supset U$ and $U$ is open.  The proof is simple and so is omitted.\\

Clearly if $A$ is a $ \wedge_\tau$-set  then $A$ is $ g\wedge_\tau$-set. But the converse is not always true  as shown in the following example.\\

\textbf{Example 3.17:}  Let $X = R$,    $ G_i $ be the all countable subsets of irrationals in $R$ and $ \tau $ = $\{X,  \emptyset, G_i\}$. Therefore $ (X, \tau) $ is a space but not a topological space. Suppose $ A = R - Q $. Therefore $ A_\tau^\wedge = X$ and $\overline{A} = X $.  Therefore $A$ is $ g\wedge_\tau$-set, since the only closed set containing $A$ is $X$. But $ A = R - Q \not = X = A_\tau^\wedge $. This implies that A is not $\wedge_\tau$-set.\\

\textbf{Lemma 3.18} :  Let $A , B$ be subsets of $X$.  Then the  following hold:

(1)$\emptyset_\tau^\wedge = \emptyset $, \qquad      $\emptyset_\tau^\vee = \emptyset$, \qquad$X_\tau^\wedge = X$,\qquad$X_\tau^\vee = X$.

(2) $A \subset A_\tau^\wedge$,   \qquad $A_\tau^\vee \subset A$.

(3)  $(A_\tau^\wedge)^\wedge = A_\tau^\wedge$,\qquad  $(A_\tau^\vee)^\vee = A_\tau^\vee$.

(4)  $A \subset B \Rightarrow A_\tau^\wedge \subset B_\tau^\wedge$.

(5)  $A \subset B \Rightarrow A_\tau^\vee \subset B_\tau^\vee$.

(6)    $(X - A)_\tau^\wedge = X - A_\tau^\vee$,  \qquad $(X - A)_\tau^\vee = X - A_\tau^\wedge$.\\

\textbf{Theorem 3.19} : Arbitrary union of $ \vee _\tau $-sets is a $ \vee_\tau $-set.

\begin{proof}
 Suppose that $ A_i $'s are $ \vee_\tau $-sets, $ i\in I $ where $ i $ is an index set and $ A=\cup\{A_i: i\in I \} $. So $ A_i\subset A $ for each i. Therefore   $  Ai_\tau^\vee \subset A_\tau^\vee $ for all $ i\in I $ and so $ \cup \{Ai_\tau^\vee: i\in I\}\subset A_\tau^\vee $. So, $ \cup \{Ai: i \in I\}= \cup\{Ai_\tau^\vee: i\in I\}\subset A_\tau^\vee  \subset A $ by Lemma 3.18 (2). Therefore $ A=A_\tau^\vee $.
 \end{proof}

\textbf{Corollary 3.20} : Arbitrary intersection of $ \wedge_\tau $-sets is a $ \wedge_\tau $-set. 
\begin{proof}
Let $\{B_i :i\in \Lambda\}$ be any collection of $ \wedge_\tau $-sets and let $ B = \cap_{i\in \Lambda} B_i $. Then $ X - B= X - \cap_{i\in \Lambda}B_i = \cup_{i\in \Lambda}(X - B_i) $. Since each $ X - B_i $ is $ \vee_\tau $-sets, by Theorem 3.19 $ X - B $ is $ \vee_\tau $-sets and so $ B $ is $ \wedge_\tau $-set.
\end{proof}

\textbf{Theorem 3.21} : Intersection of two  $ \vee _\tau $-sets is a $ \vee_\tau $-set.
\begin{proof}
 Let $A , B$ be two $ \vee_\tau $-sets. Then $ A=A_\tau^\vee, B=B_\tau^\vee $.  Now $ A\cap B \subset A $ and $ A\cap B \subset B$. Therefore $(A\cap B)_\tau^\vee\subset A_\tau^\vee $ and $(A\cap B)_\tau^\vee\subset B_\tau^\vee $. So $  (A\cap B)_\tau^\vee \subset A_\tau^\vee\cap B_\tau^\vee=A\cap B $.
 
  Conversely, suppose that $ x\in A_\tau^\vee\cap B_\tau^\vee $. Since $ x\in A_\tau^\vee $ and $x \in B_\tau^\vee$, there exist closed sets $F , P$ such that $ x\in F\subset A, x\in P \subset B $. Therefore  $ x\in F\cap P \subset A\cap B $. This implies that $x\in (A\cap B)_\tau^\vee$.  Therefore $ A_\tau^\vee \cap B_\tau^\vee \subset(A\cap B)_\tau^\vee $. So $ A\cap B= A_\tau^\vee \cap B_\tau^\vee =(A\cap B)_\tau^\vee$.
\end{proof}

\textbf{Corollary 3.22} :  Let $(\ X, \tau)$ be a space. Then the collection of all $\vee_\tau$-sets forms a topology. 
The proof follows from above two Theorems 3.19 and 3.21 and Lemma 3.18(1).\\

\textbf{Remark 3.23 } (i) : Note that in view of Theorem 3.3 a set $A$ is $ g^*$-closed if and only if there is a closed set $F$ containing $A$ such that $F\subset A_\tau^\wedge$.

(ii) : Clearly if $A$ is a $\wedge _\tau $-set,  then $A$ is $g^*$-closed if and only if $A$ is closed.  In particular $(A_\tau^\wedge)^\wedge = A_\tau^\wedge$ and so $A_\tau^\wedge$ is a $\wedge_\tau$-set.  Therefore $A_\tau^\wedge$ is $g^*$-closed if and only if $A_\tau^\wedge$ is closed. \\

\textbf{Theorem 3.24} :  If $A_\tau^\wedge$ is $g^*$-closed,  then $A$ is $g^*$-closed.
\begin{proof}
Let $A \subset X$ and $A_\tau^\wedge$ be $g^*$-closed. Then by Remark  3.23(i), there exists a closed set $F$ containing $A_\tau^\wedge$ i.e. $F \supset A_\tau^\wedge \supset A$ such that $F  \subset (A_\tau^\wedge)^\wedge = A_\tau^\wedge$, by Lemma 3.18 (3). So by Remark 3.23(i), $A$ is $g^*$-closed.
\end{proof}

\textbf{Corollary 3.25} :  If $ A_\tau^\vee $  is $ g^*$-open ,    $A$ is $ g^*$-open. 
\begin{proof}
If $ A_\tau^\vee $  is $ g^*$-open then $ X - A_\tau^\vee $ is $g^*$-closed i.e. $ (X - A)_\tau^\wedge $ is $g^*$-closed. So by above Theorem 3.24, $ X - A $ is  $g^*$-closed and hence $A$  is $ g^*$-open.    
\end{proof}

But the converse of the above theorem and corollary may  not be true in a space as shown in the following example.\\

\textbf{Example 3.26} : Suppose that $X = R - Q$,  $ G_i $'s are the countable subsets of $ X $, each contains $ \sqrt{2}$,  and $\tau = \{X, \emptyset,  G_i, A_i\} $ where $ A_i$'s are the cocountable subsets of  $ X $ each contains $ \sqrt{2} $. Then $(X, \tau)$ is a space but not a topological space.  Let $A$ be a countable infinite subset of $X$ excluding $\sqrt{2}$.  Then $A$ is a closed set and so $A$ is $g^*$-closed.  But $A_\tau^\wedge=\cap \{ \{\sqrt{2}\}\cup(X -\{\alpha\}),\alpha\in X - A,\alpha \not=\sqrt{2}\}  =\{ \sqrt{2}\}\cup A $ is  an open set which is not closed, since $ \sqrt{2}\in  A_\tau^\wedge $.  Since $ A_\tau^\wedge $ is a $ \wedge_\tau $-set, by Remark 3.23(ii), $A_\tau^\wedge$ is not $g^*$-closed. Therefore $ X-A $ is $ g^* $-open, but $ (X-A)_\tau^\vee=X-A_\tau^\wedge $ is not $ g^* $-open.\\

\textbf{Lemma 3.27 :} If $A , B$ are two subsets of $X$, then $ (A\cup B)_\tau^\wedge= A_\tau^\wedge\cup B_\tau^\wedge$.

\begin{proof}
Let $A , B$ be two subsets of $X$. Then  $ A\subset A\cup B $ implies that $ A_\tau^\wedge\subset (A\cup B)_\tau^\wedge $ and $ B\subset A\cup B $ implies that $ B_\tau^\wedge\subset (A\cup B)_\tau^\wedge $. Therefore $ (A_\tau^\wedge\cup B_\tau^\wedge)\subset (A\cup B)_\tau^\wedge $. Again, $ A_\tau^\wedge=\cap \{U_i: U_i\supset A,  U_i $ is open\} and $ B_\tau^\wedge=\cap \{V_j:  V_j\supset B,  V_j$ is open\}.  Therefore $ A_\tau^\wedge\cup B_\tau^\wedge = \cap\{(U_i\cup V_j): U_i \supset A, V_j\supset B, U_i, V_j $ are open\} $\supset\cap \{G: G\supset A\cup B, G$ is open \}=$ (A\cup B)_\tau^\wedge $. Therefore $ A_\tau^\wedge\cup B_\tau^\wedge=(A\cup B)_\tau^\wedge $.
\end{proof}

\textbf{ Theorem 3.28 :} Union of finite number of $ \wedge_\tau $-sets is a $ \wedge _\tau $-set.

The proof is straightforward by Lemma 3.27 and so is omitted.

But arbitrary union of $ \wedge_\tau $-sets is not always a $ \wedge_\tau $-set as can be seen from the example given below.\\

\textbf{ Example 3.29 :} Let $ X = R-Q $,    $ G_i $ be the all countable subsets of $X$ and $ \tau $ = $\{X,  \emptyset, G_i\}$. Then $ (X, \tau) $ is a space but not a topological space. Every singleton is an open set, so a $ \wedge_\tau $-set. Now the set $ A=[1, 2]- Q $ is not open and $ A_\tau^\wedge=X\not=A $. So $ A $ is not a $ \wedge_\tau $-set. But $ A=\cup \{\{r\}: r\in A\} $ where $ \{r\} $ is $ \wedge_\tau $-set.\\

\textbf{Theorem 3.30 :} Union of two $ g\wedge_\tau $-sets is $ g\wedge_\tau $-set.
\begin{proof}
Suppose $A , B$ be two $ g\wedge_\tau $-sets of X,  then $ A_\tau^\wedge\subset F_1 $ whenever $F_1$ is closed and $ A\subset F_1 $  and $ B_\tau^\wedge\subset F_2 $ whenever $F_2$ is closed and $ B\subset F_2 $. Let $ (A\cup B)\subset F , F$ is closed, then $ A\subset F $ and $ B\subset F $ which imply that $ A_\tau^\wedge \subset F $ and $ B_\tau^\wedge\subset F $. So $( A_\tau^\wedge\cup B_\tau^\wedge)\subset F$.  Therefore by Lemma 3.27, $(A\cup B)_\tau^\wedge = (A_\tau^\wedge\cup B_\tau^\wedge)\subset F$.  Hence $ A\cup B $ is $ g\wedge_\tau $-set.  
\end{proof}

\textbf{Corollary  3.31:} Intersection of two $ g\vee_\tau $-sets is $ g\vee_\tau $-set.
\begin{proof}
Let $A , B$ be two $ g\vee_\tau $-sets. Therefore $ (X-A) $ and $(X-B)$ are two $ g\wedge_\tau $-sets.  Then by Theorem 3.30, $ (X-A)\cup (X-B)$ is a $ g\wedge_\tau $-set i.e. $ X-(A\cap B) $ is a $ g\wedge_\tau $-set. Hence $ A\cap B $ is a $ g\vee_\tau $-set.
\end{proof}

But union of two $ g\vee_\tau $-sets is not always $ g\vee_\tau $-set as shown by the example given below.\\

\textbf{Example 3.32 :} 
Let $ X = R-Q $,    $ G_i $'s be the nonempty  countable subsets of $X - \{\sqrt{3}\}$  and let $ \tau $ = $\{X,  \emptyset, G_i\cup \{\sqrt{3}\}\}$. Then $ (X, \tau) $ is a space but not a topological space. Assume $ A=\{\sqrt{3}\} $, then $ A_\tau^\vee = \emptyset $. Only open set contained in $ A $ is $ \emptyset $ and $ \emptyset\subset \emptyset =A_\tau^\vee $. So $ A $ is a $ g\vee_\tau $-set. Again suppose $ B= \{\sqrt{5}\} $, then $ B_\tau^\vee = \emptyset$. The only open set contained in $ B $ is the set $ \emptyset $ and $ \emptyset \subset B_\tau^\vee $. So $ B $ is a $ g\vee_\tau $-set. Now suppose $ C=A\cup B=\{\sqrt{5}, \sqrt{3}\} $. Then $ C_\tau^\vee =\emptyset $. Also $ C $ is an open set and $ C\subset C $ but $ C\not\subset C_\tau^\vee  $. Therefore $ A\cup B $ is not a $ g\vee_\tau $set.\\

Similarly, intersection of two $ g\wedge_\tau $-sets is not always a $ g\wedge_\tau $-set as revealed by the undernoted example.\\

\textbf{Example 3.33 :} Consider the space $ (X,  \tau )$ as in Example 3.32 where $ A=\{\sqrt{3}\} $ and $ B=\{\sqrt{5}\} $ are $ g\vee_\tau $-sets and $ A\cup B $ is not a $ g\vee_\tau $-set. Therefore $ X-(A\cup B) $ is not a $ g\wedge_\tau $-set. Now $ X-A $ and $ X-B $ are $ g\wedge_\tau $-sets. So $ (X-A)\cap (X-B)=X-(A\cup B) $ is not a $ g\wedge_\tau $-set.\\

\textbf{Remark 3.34} : Clearly if $A$ is $ \vee_\tau$-set ,  then $A$ is $ g\vee_\tau$-set. But the converse is not necessarily true  as shown in the following example. \\

\textbf{Example 3.35} : Let $X = R$,    $ G_i $'s be the countable subsets of irrationals in $R$ and let $ \tau $ = $\{X,  \emptyset, G_i\}$. Then $ (X, \tau) $ is a space but not a topological space. Suppose $A = Q$.  Then $ (X - A)_\tau^\wedge = X $ and            $ \overline{X - A} $ = $ X $.  Therefore $ X - A$ is a $ g\wedge_\tau$-set, since the only closed set containing $ X - A$ is X.  So $A$ is $ g\vee_\tau$-set . Now $ X - A$ = $ R - Q \not= X = \ (X - A)_\tau^\wedge $. So $ X - A$ is not $ \wedge _\tau$-set.  Therefore $A$ is not $ \vee_\tau $-set.\\

\textbf{Theorem 3.36 :}   For each $ x\in X $,  $ \{x\} $ is either open or a $ g\vee_\tau $-set.

The proof is straightforward and so is omitted. \\
                    
\textbf{Note 3.37} :  In view of above theorem it follows that for any $x \in X$, either $\{x\}$ is $ g^* $-open  or  $\{x\}$  is $ g\vee_\tau $-set.\\

\textbf{Theorem 3.38}  :  Let $A$ be a $ g\vee_\tau $-set  in a space $ (X,\tau)$. If $ A_\tau^\vee \cup(X - A) \subset F, F$ is closed,  then $F = X$.
\begin{proof}
 Let $ A $ be a $ g\vee_\tau $-set in a space $ (X, \tau) $. Then $ U\subset A_\tau^\vee $ whenever $ U\subset A $ and $ U $ is open. Let $ F  $ be a closed set such that $ A_\tau^\vee \cup (X - A)\subset F $. This implies $ F^c \subset(A_\tau^\vee \cup A^c)^c = (A_\tau^\vee)^c \cap A $ i.e. $ X - F \subset (X - A_\tau^\vee)\cap A $......(1).  Therefore $ X - F \subset A $. Since $ X - F $ is open and $ A $ is $ g\vee_\tau $-set, then $ X - F \subset A_\tau^\vee $. From (1) we have $ X - F \subset X - A_\tau^\vee= (A_\tau^\vee)^c $. Therefore $ X - F \subset A_\tau^\vee \cap (A_\tau^\vee)^c = \emptyset $ which implies $ X = F $.
\end{proof}

\textbf{Corollary 3.39}  : Let $A$ be a $g\vee_\tau $-set in a space $(X, \tau)$.  Then $A_\tau^\vee \cup(X - A)$ is closed if and only if $A$ is a $ \vee_ \tau $-set .  
\begin{proof}
Let $ A $ be a $ g\vee_\tau $-set in a space $ (X, \tau) $ and let $A_\tau^\vee \cup(X-A)$ be closed. Then by Theorem 3.38, $A_\tau^\vee \cup(X - A)= X $ which implies that $(A_\tau^\vee)^c \cap A = \emptyset$. So $ A \subset \{(A_\tau^\vee)^c\}^c = A_\tau^\vee $.  Again $ A_\tau^\vee \subset A $.  Therefore $ A = A_\tau^\vee$ which  implies that $ A $ is $ \vee_\tau $-set.

Conversely, Let $ A $ be a $ g\vee_\tau $-set and also $ \vee_\tau $-set.  Therefore $ A \subset A_\tau^\vee $, since $A = A_\tau^\vee$. So $ (A_\tau^\vee)^c \cap A = \emptyset $ which implies that $ A_\tau^\vee \cup (X - A) = X $. Therefore  $A_\tau^\vee \cup(X - A)$ is closed.
\end{proof}

\textbf{Theorem  3.40} :  If $A$ is a $ \vee_\tau $-set, then $A$ is $ g^*$-open  if and only if $A$ is open. 
\begin{proof}
Let $ A $ be open. Then clearly $ A $ is $ g^* $-open. Conversely, Let $ A $ be a $ \vee_\tau $-set and let $ A $ be $ g^* $-open. So by Theorem 3.14 there is an open set $ U\subset A $ such that $ A_\tau^\vee \subset U $. This implies that $ A\subset U, $ since $ A = A_\tau^\vee $. Therefore $ A = U $ and so $ A $ is open.
\end{proof}
In particular, since $(A_\tau^\vee)^\vee = A_\tau^\vee$, so $A_\tau^\vee$ is a $\vee_\tau $-set.  Therefore $A_\tau^\vee$ is $g^*$-open if and only if $A_\tau^\vee$ is open. \\

\section{\bf $T_\omega $-Space}

\textbf{Definition 4.1}  \cite{DP} :  A space $ (X, \tau) $ is said to be $ T_\omega $-space if every $ g^*$-closed  set is closed.\\

  In Theorem  16  \cite{DP}, it is shown that every $ T_\omega $-space is $ T_0 $-space. But the converse is not true as shown in the Example 6 \cite {DP}. Also, in Examples 6 and 7 \cite{DP} it has been shown  that $ T_\omega $ and $ T_1 $ axioms in a space are independent of each other.\\

\textbf{Definition 4.2}  \cite{DP}) : For any $ E \subset X $  let $ \overline{E^*}=\cap\{A:E\subset A , A $ is  $ g^*$-closed  set in $ X $\},  then $\overline{E^*}$ is called $ g^*$-closure  of $E$.\\

We consider the following sets which will be used frequently in the sequel:

$C =\{A: \overline{(X-A)}$ is closed\} and $ C^*=\{A:\overline{(X-A)^*}$ is $ g^*$-closed\}.\\

\textbf{Theorem 4.3}  \cite{DP}  ):  A space $ (X, \tau) $ is $ T_\omega $ if and only if 

(a)  for each $ x \in X $,  $ \{x\} $ is either open or closed and

(b)  $ C = C^* $.   \\

\textbf{Theorem 4.4} : The following are equivalent:
 
(1)  $ (X,  \tau) $ is $ T_\omega$-space.

(2)  Every $ g\wedge_\tau $-set  is $ \wedge_\tau $-set  and $ C=C^* $. 

(3)  Every $ g\vee_\tau $-set  is $ \vee_\tau $-set  and $ C=C^* $.
 
\begin{proof}
(1) $\Rightarrow$ (2): Let $ (X,  \tau) $ be $ T_\omega $-space and let $A$ be a $ g\wedge_\tau $-set . We wish to prove $ A \supset A_\tau^\wedge $.  If not, suppose that $ x\in A_\tau^\wedge $ but $ x\not\in A $. By Theorem 4.3, $\{x\}$ is either open or closed. We discuss two cases: 

case (i): Suppose $ \{x\} $ is open. So $ X - \{x\}$ is a closed set containing $A$. Since $A$ is $ g\wedge_\tau $-set, we get $ A_\tau^\wedge \subset (X - \{x\})$  which implies that $ x\not\in A_\tau^\wedge  $,  a contradiction. 

Case (ii): Suppose $ \{x\} $ is closed.  Then $ X -\{x\}$ is an open set containing $A$.  But $ x\in A_\tau^\wedge $=$\cap\{U,U\in \tau,U\supset A\}\subset X -\{x\}$, a contradiction. Hence  in any case $ A \supset A_\tau^\wedge$ 
and so  $ A=A_\tau^\wedge$ which implies that $ A $ is $ \wedge_\tau $-set . Also by Theorem 4.3, $  C=C^* $. 

(2) $\Rightarrow$ (3):  Let $A$ be a $ g\vee_\tau $-set  and $ C=C^* $. Then by Definition 3.15,  $ X - A $ is $ g\wedge_\tau $-set .  By supposition $ X - A $ is a $ \wedge_ \tau $-set . So $A$ is a $ \vee_\tau $-set  . 

(3) $\Rightarrow$ (2):  Let  $B$ be a $ g\wedge_\tau $-set and $ C = C^* $.  So  $ X - B $ is $ g\vee_\tau $-set.  By supposition, $ X - B $ is $ \vee_\tau $-set. So   $B$ is $ \wedge_\tau $-set.

(2) $\Rightarrow$ (1) : Let $ x \in X $ and $ C = C^* $.  We will prove $ \{x\} $ is either open or closed.  If $ \{x\} $ is not open,  $ X - \{x\} $ is not closed. So $X$ is the only closed set containing $ X - \{x\} $.   Also   $ (X - \{x\})_\tau^\wedge \subset X $ .Therefore by definition,  $ X - \{x\} $ is a $ g\wedge_\tau $-set . By supposition,  $ X - \{x\} $ is $ \wedge_\tau $-set  i.e. $ (X - \{x\})_\tau^\wedge $  = $ X - \{x\} $.  So $ X - \{x\} $ must be open. So $\{x\}$ is a closed set. Therefore the space is $ T_\omega $.
\end{proof}

\section{\bf $ \lambda^*$-closed sets and $ \lambda^* $-open sets in a space}

\textbf{Definition 5.1 }:  A subset $A$ of a space $ (X,  \tau) $  is said to be $ \lambda^*$-closed   if $ A=L\cap \overline{F} $  where $L$ is a $ \wedge_\tau $-set   and $F$ is a subset of $X$.

Clearly every $ \wedge _\tau $-set  is $ \lambda^*$-closed and closed set is $ \lambda^*$-closed. But the converse may not be
true as shown in the following example.\\

\textbf{Example 5.2 :} Let $X=R-Q$,  $ \tau = \{X,  \emptyset,  G_i\}$,  where $ G_i $  be the all countable subsets of $X$. So $ (X,  \tau) $ is a space but not a topological space.  Let $A$ be the set of all irrational numbers in $ (0,  \infty) $. Then clearly $A$ is $ g^*$-closed but $ A $ is not closed. Now $ A_\tau^\wedge=X $  and $\overline{A}=A $. Hence $ A_\tau^\wedge \cap \overline{A}=A$.  So $A$ is $ \lambda^*$-closed . But $ A_\tau^\wedge\not=A$  which implies $A$ is not $ \wedge_\tau $-set.\\

\textbf{ Lemma  5.3} :  For a subset
 $A$ of a space  $ (X,  \tau) $  
the following are equivalent: 

(i)  $A$ is $ \lambda^*$-closed.

 (ii) $ A=A_\tau^\wedge\cap \overline{F},F\subset X$.

(iii)  $ A=L\cap \overline{A}$,  L is a $ \wedge_\tau $-set . 

(iv)  $ A=A_\tau^\wedge\cap \overline{A} $. 

\begin{proof}
(i) $ \Leftrightarrow$ (ii): Let $A$ be $ \lambda^* $-closed i.e. $ A=L\cap \overline{F} $ where $ L=L_\tau^\wedge, F\subset X $. Then $ A\subset L $ which implies that $ A_\tau^\wedge \subset L_\tau^\wedge $ and $ A\subset \overline{F} $. Therefore $ A\subset A_\tau^\wedge \cap \overline{F}\subset L_\tau^\wedge \cap \overline{F}=L\cap \overline{F} =A $. Hence  $ A= A_\tau^\wedge\cap \overline{F} $.

 Conversely, suppose $ A_\tau^\wedge\cap \overline{F}= A $ where $ F\subset X $. Since $ A_\tau^\wedge $ is a $ \wedge_\tau $-set then    $ A  $ is $ \lambda^*$-closed.

(i) $ \Leftrightarrow $ (iii):  Assume $ A=L\cap \overline{F} $ where $ L=L_\tau^\wedge $ and $ F\subset X $. Then $ A\subset L $ and $ A\subset \overline{F}$. So $ \overline{A}\subset \overline{\overline{F}}=\overline{F} $. Therefore $ A\subset L\cap\overline{A}\subset L\cap \overline{F}=A $ and hence $ A= L\cap \overline{A} $. 

Converse part directly follows from definition.

(i) $ \Leftrightarrow (iv) $: Assume $ A=L\cap \overline{F} $ where $ L=L_\tau^\wedge $ and $ F\subset X $.  Then $ A\subset L $, so $ A_\tau^\wedge \subset L_\tau^\wedge$. Again $ A\subset \overline{F} $. This implies that $ \overline{A}\subset \overline{\overline{F}}=\overline{F}$.  Therefore $ A\subset A_\tau^\wedge \cap\overline{A}\subset L_\tau^\wedge \cap\overline{F}=L\cap \overline{F} = A $ and hence $ A=A_\tau^\wedge\cap \overline{A} $.

 Converse part directly follows from definition.
\end{proof}

\textbf{Theorem 5.4} : If  $A$ is $ \lambda^*$-closed  and $ g\wedge_\tau $-set such that $ \overline{A} $ is closed,  then $A$ is a $ \wedge_\tau $- set.
\begin{proof}
Suppose that $A$ is $ \lambda^*$-closed and $ g\wedge_\tau $-set such that $ \overline{A} $ is closed. Since $ \overline{A} $  contains $ A $ and  $ A $ is $ g\wedge_\tau $-set, then $ A_\tau^\wedge \subset \overline{A} $. Since $ A $ is $ \lambda^* $-closed,  $ A = A_\tau^\wedge\cap \overline{A} $ which implies $ A = A_\tau^\wedge $. Hence $ A $ is $ \wedge_\tau $-set.
\end{proof}

\textbf{Remark 5.5} : From Lemma 5.3 (iv)  we can say that a subset $A$ is said to be $ \lambda^*$closed   if $A$ can be expressed as the intersection of all open sets and all closed sets containing it.

In view of above remark 5.5 and by Theorem 18 \cite{DP} we have the following theorem 5.6.  However we are also giving  a separate proof of the Theorem 5.6.\\

\textbf{Theorem 5.6} : A space $ (X,  \tau) $ is $ T_\omega $ if and only if  every subset of $(X,\tau)$ is $\lambda^*$-closed and $C = C^*$.

\begin{proof} 
 Suppose every subset of  $(X,\tau)$ is $\lambda^*$-closed and $C = C^*$ and $x$ $\in$ X. We  shall show that $\{x\}$ is either open or closed. Suppose $\{x\}$ is not open, then $X-\{x\}$ is not closed. Since $X-\{x\}$ is also a $\lambda^*$-closed set then $X-\{x\} = (X-\{x\})_\tau^\wedge\cap(\overline{X-\{x\}})$, (by Lemma 5.3(iv)) = $(X-\{x\})_\tau^\wedge\cap X =(X- \{x\})_\tau^\wedge$. Therefore $X-\{x\}$ is a $\wedge_\tau $-set. So $X-\{x\}$ is an open set  which implies $\{x\}$ is closed. Then by Theorem 4.3, $(X,\tau)$ is $T_\omega$-space.
 
 Conversely, suppose that $ (X, \tau) $ is $T_\omega$-space and $ A\subset X $.  Then by Theorem 4.3, every singleton is either open or closed and  $C = C^*$.  So each $ x\in X - A $,  either $ \{x\}\in \tau $  or  $ (X - \{x\}) \in \tau $. Let $ A_1=\{x : x\in X - A, \{x\}\in \tau \} $,  $ A_2=\{x : x\in X - A, X - \{x\}\in\tau \}$,  $ L=\cap[X - \{x\} : x\in A_2]=X - A_2 $ and   $F=\cap[X - \{x\} : x\in A_1]=X - A_1 $. 
 
 Note that L is a $ \wedge_\tau $-set   i.e. $ L=L_\tau^\wedge$  and $ F=\overline{F} $.
 
 Now,  $L\cap\overline{F}=(X - A_2)\cap (X - A_1)=X - (A_1\cup A_2)=X -(X - A)=A$.  Thus $A$ is  $\lambda^*$- closed. 
\end{proof}

\textbf{Remark 5.7 :}  It is already seen that if a subset $A$  is closed then $A$ is $ g^*$-closed  and $ \lambda^*$-closed. But the converse of this result may not hold in a space as seen in Example 5.2,  although it is true in a generalised topological space \cite{MS}. However,  it is true in a space if the additional condition that $ C=C^* $  holds  which is shown  in the following Theorem 5.8. \\
 
\textbf{Theorem  5.8 :} If $A$ is $ g^*$-closed  and $ \lambda^*$-closed  and satisfies the condition $ C=C^* $, then $A$ is a closed set.

\begin{proof}
Suppose $A$ is $ g^*$-closed  and $ \lambda^*$-closed  satisfying the condition $ C=C^* $. Since $ A $ is $ g^* $-closed,  there is a closed set $F$ containing $A$ such that $ F\subset  A_\tau^\wedge $.  Since $ A\subset F,  \overline{A}\subset\overline{F}$ and so $ \overline{A}\subset F\subset A_\tau^\wedge $.  Again since $A$ is $ \lambda^*$-closed,  $A=A_\tau^\wedge \cap\overline{A}=\overline{A}$...............(1). 

 Now $A$ is $ g^*$-closed, so $ \overline{A^*}=A$, therefore $ (X - A) \in C^* $.  Since $ C=C^* $, $ (X - A) \in C $ which implies that $\overline{A} $ is closed. Therefore, by (1),   $A$ is a closed set.  
\end{proof}

\textbf{Remark  5.9 :}  Union of two $ \lambda^*$-closed  sets may not be $ \lambda^*$-closed  set as revealed in the following example.\\

\textbf{Example  5.10 :}  Let $X=R-Q$,  $ \tau = \{X,  \emptyset,  \{\sqrt{3}\}, G_i\cup \{\sqrt{3}\}\}  $ where $ G_i $  be  all countable subsets of $X-\{\sqrt{2}\}$. Then $ (X,  \tau) $ is a space but not a topological space. Suppose $A =\{\sqrt{2}\} $,  then  $A_\tau^\wedge$ = X and $\overline{A}=A$ and   so $ A=A_\tau^\wedge\cap\overline{A}$. This implies that $ A $ is $ \lambda^*$-closed. Again suppose that $B=\{\sqrt{3}\}$, then B is also $ \lambda^*$-closed, since $ B = B_\tau^\wedge $ and $ \overline{B} = X $. Now  let $ C=\{\sqrt{2},\sqrt{3}\} $. Then $ C_\tau^\wedge=X $ and $ \overline{C}=X $. So $ C_\tau^\wedge\cap\overline{C}=X\not=C $ which implies that $ C $  is not $ \lambda^*$-closed.\\

\textbf{Remark  5.11} : It is shown in \cite{DP} that in a bicompact space (Alexandroff space), $ g^* $-closed sets may not be bicompact and likewise,  $ \lambda^* $-closed sets in a bicompact space may not be bicompact as shown in the following example.\\

\textbf{Example 5.12} : Let $X=R-Q$,  where $ G_i $  be the all countable subsets of $X-\{\sqrt{2}\}, A_i $ be the all cofinite subsets of $ X $ and $ \tau = \{X,  \emptyset,   G_i, A_i \} $. Then $ (X,  \tau) $ is a space but not a topological space. Clearly $(X, \tau)$ is a bicompact space. Suppose $ A $ is the set of all irrationals in $ (0, 1)$. Then there is no open set of $ G_i $ type containing $ A $ and so  $  A_\tau^\wedge =\cap \{X-\{\alpha\}: \alpha\in X- A\} = A $. Therefore $ A $ is a $ \wedge_\tau $-set  which implies $ A $ is a $ \lambda^* $-closed set. But $ A $ is not bicompact since $ \{\{r\} :r\in A\} $ forms an open cover for $ A $ which has no finite subcover.\\

\textbf{Definition 5.13 }:  A subset $ A $ of $ X $  is said to be $\lambda^*$-open  set if $ X - A $ is $ \lambda^*$-closed set.\\

\textbf{Theorem 5.14} : A subset $ A $ of $ X $   is  $ \lambda^* $-open set if and only if  $ A=M\cup Int(V) $ where $M$ is a $ \vee_\tau $-set  and $V$ is a subset of $X$.
\begin{proof}
Let $A$ be a $ \lambda ^* $-open set. Then $ X - A $ is a $ \lambda^* $-closed set. So $ X - A =L\cap\overline{F} $, where $L$ is a $ \wedge_\tau $-set and $ F\subset X $. Therefore  $ X-L\subset A $ and  $ X-\overline{F}\subset A $. So by Theorem 2.8, $ A = L^c\cup (\overline{F})^c = M\cup Int(X-F) = M\cup Int(V)$ where $ M = X-L $, a $ \vee_\tau $-set and $ V = X-F \subset X $.  

Conversely,  Let $ A=M\cup Int(V) ,  M $ is a $ \vee_\tau $-set and $ V\subset X $. So $ A^c=M^c\cap (Int(V))^c $ where $ M^c $ is a $ \wedge_\tau $-set and $ (Int(V))^c   =\overline{X-V}$ by Theorem 2.8.  Therefore $ X - A=L\cap \overline{F} $, where $ L=M^c $ and $ F=X-V $.  Therefore $ X - A $ is a $ \lambda^* $-closed set and hence $A$ is a $ \lambda ^* $-open set.
\end{proof}

\textbf{Remark 5.15} : Clearly $ \vee_\tau $-sets are $ \lambda ^* $-open sets and open sets are $ \lambda^* $-open sets. On the other hand  if a set $ A $ is open and $ g\vee_\tau $-set  then  $A$ is a $ \vee_\tau $-set. Again if $A$ is $ \lambda^* $-open, $g^*$-open and satisfies the condition $ C=C^*$,  then $ A $ is open.\\

\textbf{Theorem 5.16 :} $A$ is $ \lambda^*$-open  if and only if $  A=A_\tau^\vee \cup Int(A)$.
\begin{proof}
Let $A$ be $ \lambda^* $-open. Then $ A=M\cup Int(V)$ where $M$ is a $ \vee_\tau $-set and $ V\subset X $. Since $ M\subset A$,  $ M_\tau^\vee \subset A_\tau^\vee $ and since $ Int(V)\subset A,  Int(Int(V))\subset Int(A) $. So, $ A=M\cup Int(V)=M_\tau^\vee\cup Int(V)\subset A_\tau^\vee \cup Int(V)=A_\tau^\vee\cup Int(Int(V))\subset A_\tau^\vee\cup Int(A) $. Again since  $ A_\tau^\vee\subset A $ and $ Int(A) \subset A, A_\tau^\vee\cup Int(A)\subset A $. Therefore we get $ A=A_\tau^\vee\cup Int(A) $.

Conversely,  suppose  $ A=A_\tau^\vee\cup Int(A) $. Since $ A_\tau^\vee $ is a $ \vee_\tau $-set and $ A\subset X $, $A$ is $ \lambda^* $-open, by Theorem 5.14.
\end{proof}

 In  \cite{MS} it is seen that the collection of all $ \lambda_\mu $-open sets forms a generalised topology $ \mu $ on $ X $, but unlikely the collecion of $ \lambda^* $-open sets does not form a space structure $ \sigma $ on $ X $  as shown in the following example.\\

\textbf{Example 5.17} : Intersection of two $ \lambda^* $-open sets may not be a $ \lambda^* $-open set.

Consider the sets $ X, A, B, C $ and $ \tau $ as in Example 5.10. Then $A ,  B $ are $ \lambda^* $ -closed sets and so $ X - A $ and $ X - B $ are $ \lambda^* $-open sets. But $ C=A\cup B $ is not  $ \lambda^* $-closed. Thus $ (X-A)\cap (X-B)=X-(A\cup B) $ is not $ \lambda^* $-open. 

\section{\bf $ \textit T_\frac{\omega}{4}$-space, $\quad \textit T_\frac{3\omega}{8}$-space,   $ \quad\textit  T_\frac{5\omega}{8}$-space  }

\textbf{Definition 6.1} :  A space $ (X,  \tau) $  is called $T_\frac{\omega}{4}$-space  if for every finite subset $P$ of $X$ and for every $ y\in X - P $,  there exists a set $A_y$ containing $P$ and disjoint from $\{y\}$  such that $A_y$ is either open or closed.\\

\textbf{Theorem 6.2} :    Every $T_\frac{\omega}{4}$-space  is a $T_0$-space.
\begin{proof}
Suppose that  $ x , y $ are two distinct points in $ X $.  Since the space is $T_\frac{\omega}{4}$, then for every $ y\in X -\{x\} $ there exists a set $A_y$ such that for every  $\{x\}\subset A_y$  and $ \{y\} \not \subset A_y $ where $ A_y$ is either open or closed. This implies  $ T_\frac{\omega}{4}$-space  is  $T_0$-space.  
\end{proof}
But the converse may not be true as shown in the Example 6.5.\\

\textbf{Theorem  6.3} :  A space $ (X,  \tau) $ is $ T_\frac{\omega}{4}$ if and only if every finite subset of $X$ is $ \lambda^*$-closed.  

\begin{proof}
Suppose $ (X,  \tau) $ is $ T_\frac{\omega}{4}$-space  and $P$ is a finite subset of $X$.  So for every $ y\in X - P $ there is a set $ A_y $ containing $P$ and disjoint from $ \{y\} $ such that $ A_y $ is either open or closed.  Let $L$ be the intersection of all open sets $ A_y $ and $F$ be the intersection of all closed sets $A_y$.  Then $ L=L_\tau^\wedge$ i.e. $L$ is a $ \wedge_\tau $-set  and $ \overline{F}=F$. Therefore $ P=L\cap F=L\cap \overline{F}$. So $ P $ is $ \lambda^*$-closed. 

Conversely,  let $P$ be a finite set. So by the condition it is $ \lambda^* $-closed.  Then by Lemma 5.3(iii) $ P=L\cap \overline{P} $, where $L$ is a $ \wedge_\tau $-set. Let  $ y\in X - P $. If $ y\not\in \overline{P} $ then there exists a closed set $ F=A_y $ containing $P$ such that $ \{y\} \not\subset A_y $.  Again if $ y\in \overline{P} - P $, then $ y\not\in L$  and so $ y\not\in U $ for some open set $ U = A_y $ containing $L$, since $ L_\tau^\wedge=L $.  So $ P\subset U $.  Hence $ (X,  \tau) $ is $ T_\frac{\omega}{4}$-space.
\end{proof}

\textbf{Theorem  6.4} : A space $ (X,  \tau) $ is $ T_0 $ if and only if every singleton of $X$ is $ \lambda^*$-closed. 

\begin{proof}
Let the  space $ (X,  \tau) $ be  $ T_0$ and $ x\in X $.  Take a point $ y \in X - \{x\} $. So there exists a set $ A_y $ containing $ \{x\} $ but $ y\not \in A_y $ such that $ A_y $ is either open or closed. Let $L$ and $F$ be the intersection of all such open sets $ A_y $ and all such closed sets $ A_y $ respectively. Then $ \overline{F}=F $ and so $ \{x\} = L\cap F=L\cap\overline{F}$ where $L$ is a $ \wedge_\tau $-set. This implies that $ \{x\}$ is a $ \lambda^*$-closed  set by definition. 

Conversely,  let $x\in X$.  So $ \{x\} $ is  $ \lambda^*$-closed. Then $ \{x\} = \{x\}_\tau^\wedge\cap \overline{\{x\}}$  by  Lemma 5.3(iv).  Suppose $ y\in X -\{x\} $. If $ y\not\in \overline{\{x\}} $,  then there exists a closed set $F$ containing $ \{x\} $ but not containing $y$. If $ y\in \overline{\{x\}} $ then $  y\not \in \{x\}_\tau^\wedge$. Thus $ y\not \in V $ for some open set $V$ containing $ \{x\} $.  Hence the space $ (X,  \tau) $ is $ T_0 $. \end{proof}

\textbf{Example  6.5} : Example of a $ T_0$-space which is not $ T_\frac{\omega}{4} $.

  Consider the space $ (X, \tau) $, the subsets $ A,B,C $ as in the Example 5.10 where we see that every singleton in $X$ is $ \lambda^*$-closed. For, let $ r \in X, r\not = \sqrt{2}, \sqrt{3} $ then $\{r\}_\tau^\wedge = \{\sqrt{3}, r\}$ and $ \overline{\{r\}} = \{\sqrt{2}, r\} $. Therefore $\{r\} = \{r\}_\tau^\wedge \cap \overline{\{r\}}$. So \{r\} is $ \lambda^*$-closed.  Therefore the space $ (X,  \tau) $ is $T_0$-space. But it is proved there that the finite subset $C$ is not a $ \lambda^*$-closed  set.  So by Theorem 6.3,  $ (X,  \tau) $  is not $ T_\frac{\omega}{4}$-space.\\

\textbf{Definition 6.6} : A space$ (X,  \tau) $  is called {$T_\frac{3\omega}{8}$-space if for every countable subset $P$ of $X$ and for every $ y\in X - P$,  there exists  a set $ A_y $ containing $P$ and disjoint from $ \{y\} $ suh that $ A_y $ is either open or closed.\\

 \textbf{Theorem  6.7} :  A space$ (X,  \tau) $ is {$T_\frac{3\omega}{8}$-space if and only if every countable subset of $X$ is $ \lambda^*$-closed.
 Proof is similar to the proof of Theorem 6.3,  so is omitted.\\

 \textbf{Definition 6.8} : A space$ (X,  \tau) $  is called {$T_\frac{5\omega}{8}$-space if for any subset $ P $ of $ X $ and for every $ y\in X - P$,  there exists  a set $ A_y $ containing $P$ and disjoint from $ \{y\} $ such that $ A_y $ is either open or closed.\\

 \textbf{Theorem  6.9  :} A space$ (X,  \tau) $ is {$T_\frac{5\omega}{8}$-space if and only if for every subset $ E $ of $ X $ is $ \lambda^*$-closed .
  Proof is similar to the proof of Theorem 6.3,  so is omitted.\\

  Note that {$T_\frac{5\omega}{8}$ axiom does not imply $ C=C^* $.  \\

 \textbf{Remark  6.10  :} It follows from Theorem 5.6,  Theorem 6.9, Theorem 6.7,  Theorem  6.3 that every $ {T_\omega}$-space is {$T_\frac{5\omega}{8}$-space and {$T_\frac{5\omega}{8}$-space is {$T_\frac{3\omega}{8}$-space and {$T_\frac{3\omega}{8}$-space is $ T_\frac{\omega}{4}$-space. 
 
 However, the converse of each implication may not be true as shown  in the undermentioned examples.\\
 
  \textbf{Example 6.11 :} Example of a  {$T_\frac{5\omega}{8}$-space which is not $ T_\omega$-space.
 
  Let $X=R-Q$,  $ \tau = \{X,  \emptyset,  G_i\}$,  where $ G_i $  be the all countable subsets of $X$. Then $ (X,  \tau) $ is a space but not a topological space.  We see that the space $ (X,  \tau) $ is {$ T_\frac{5\omega}{8}$-space,  since any finite set,  countable set  and any uncountably infinite set are also $ \lambda^*$-closed  sets.  But the set of all irrationals in $ (0, 1) $  is  $g^*$-closed but not closed which implies the space is not $ T_\omega$-space.\\

 \textbf{Example 6.12 :} Example of a  $T_\frac{3\omega}{8}$-space which is not $ T_\frac{5\omega}{8}$-space.
 
  Let $ X=R-Q $,  $ G_i $  be  the  countable subsets of $X$ containing $ \sqrt{2} $ and $ \tau=\{X,  \emptyset,  G_i\} $. So $ (X,  \tau) $ is a space but not a topological space.  Take any countable subset $A\subset  X $. Then if  $\sqrt{2}\in A$, $A$ is an open set.  Therefore $ A=A_\tau ^\wedge $ which implies that $A$ is $\lambda^*$-closed.  If $ \sqrt{2}\not\in A,  A_\tau^\wedge=\{\sqrt{2}\}\cup A$ and  $\overline{A}=A $. This implies that $ A $  is  $ \lambda^*$-closed. So $ (X,  \tau) $ is a  {$T_\frac{3\omega}{8}$-space. Now let $B$ be an uncountably infinite subset of $X$ containing the point $ \sqrt{2} $.  Therefore $ B_\tau^\wedge = X$ and $\overline{B}=X $. Therefore $ B_\tau^\wedge\cap \overline{B}\not=B $ which implies that $ B $ is not $ \lambda^*$-closed .  So $ (X,  \tau )$ is not $ T_\frac{5\omega}{8}$-space.\\

\textbf{Example  6.13  :} Example of a  $T_\frac{\omega}{4}$-space which is not $ T_\frac{3\omega}{8}$-space.

Suppose that  $ X=R-Q $ and $X^{*}=\{2\}\cup X $. Let $  \tau=\{\emptyset, X^{*},  \{2\}\cup (X - A); A\subset X $\} where $ A $'s  are the finite subsets of $ X $.  Therefore  $ (X^*,  \tau) $ is a  topological space so a space also.  Take any finite  subset $E \subset X^{*} $,  we get the following observations:

(i)  if $ 2 \in E, E_\tau^\wedge=\cap\{ (X - \{\alpha\})\cup\{2\} ,\alpha \in X - E\}=\{2\}\cup E =E $ which implies that $ E $ is a $ \wedge_\tau $-set. Therefore $  E$ is a $ \lambda^*$-closed set .

(ii) if $ 2 \not \in E ,  E$ is a closed set which implies that $E$ is $\lambda^*$-closed.

So $ (X^{*},  \tau) $ is a $T_\frac{\omega}{4}$-space.

Now suppose that $Y$ is a countably infinite subset of $X$, so $2\not\in Y $. Here closed sets are finite. Therefore $ \overline{Y}=X^{*} $ and $ Y_\tau^\wedge=\{2\}\cup Y $. Thus $ \overline{Y}\cap Y_\tau^\wedge=\{2\}\cup Y \not=Y $. Therefore $ Y $ is not $ \lambda^{*} $-closed. Hence $ (X, \tau) $ is not $ T_\frac{3\omega}{8} $-space.\\

\textbf{Theorem  6.14} :  A space $ (X,  \tau) $ is $ T_1 $ if and only if it is $ T_0 $ and $ R_0 $. 

\begin{proof}
Let  $ (X,  \tau) $ be a $ T_1 $-space.  Obviously then it is $  T_0 $.  Let $A \subset X, A $ be open,  $ x \in A $.  Since the space is $ T_1 $, for $ x,  y \in X $ and $ x\not=y $ there are open sets $U , V $ such that $ x\in U $  and $ y \not \in U $ and $ y\in V $  and $ x\not\in V $. Hence $y$ cannot be an accumulation point of $ \{x\} $. Therefore no point lying outside $ \{x\} $ can be an accumulation point of $ \{x\} $,  so $ \{x\}'\subset  \{x\} $. Hence $ \overline{\{x\}} =\{x\}\subset A $. So  $ (X,  \tau) $ is $ R_0 $-space.

Conversely, let $ (X,  \tau) $ be $ T_0 $ and $ R_0 $. So for $ x,  y \in X $ and $ x\not=y $,  either $ x\not\in \overline{\{y\}} $ or $ y\not\in \overline{\{x\}} $. Suppose that $ x\not\in \overline{\{y\}} $.  Then there exists a closed set $F$ containing $y$ such that $ x \not\in F $. Therefore $ x\in X - F$,  an open set and $ y \not \in X - F $. Since the space is also $ R_0 $, $\overline{\{x\}} \subset X - F $. So $ \overline{\{x\}}\cap F=\emptyset $ which implies that $ \overline{\{x\}}\cap \{y\}=\emptyset $. Hence $ y \not\in \overline{\{x\}} $ and so $ y $ is not a limiting point of $ \{x\} $. Therefore there exists an open set $V$ containing $y$ such that $ x \not \in V 
$. Since $ x,  y \in X $ and $ x\not=y $,  we get two open sets $ X - F$ and $V$ containing $x , y$ respectively and $ y \not \in X - F $ and $ x \not \in V $.  Thus $ (X,  \tau) $ is $ T_1 $-space.
\end{proof}

\textbf{Theorem 6.15} :  A space $ (X,  \tau) $ is $ T_1 $ if and only if   every  singleton is $ \wedge_\tau $-set.   

\begin{proof}
Let  $ (X,  \tau) $ be $ T_1$. So by Theorem 6.14, $ (X, \tau) $ is   $  T_0$ and  $ R_0 $. Since $ (X,\tau)$ is $ T_0 $,  every singleton is $ \lambda^*$-closed, by Theorem 6.4. Suppose $ x\in X $,  then $ \{x\} $ is $ \lambda^*$-closed. So $ \{x\}=\{x\}_\tau^\wedge \cap \overline{\{x\}} $,  by Lemma 5.3 (iv). We claim that $ \{x\}= \{x\}_\tau^\wedge $.  If not, there exists $ y \in  \{x\}_\tau^\wedge - \{x\}$. So $ y\not\in \overline{\{x\}} $, and hence there is a closed set $F$,  $ F \supset \{x\}$ such that $ y \not \in F $. Therefore $ y\in X - F $, an open set. Again since $ (X, \tau) $ is  $ R_0,  \overline{\{y\}} \subset X - F $. Thus $\overline{\{y\}}\cap F=\emptyset $.  Since $ x \in F, x\not\in\overline{\{y\}}$. Therefore there exists an open set $V$ containing $x$ but $ y \not \in V $.  This implies that $ y \not\in \{x\}_\tau^\wedge $,  a contradiction.  Hence $ \{x\} $ is  a $ \wedge_\tau  $-set.

Conversely, let $ x,y \in X $ and $ x\not=y $.   So $ y \not \in \{x\} $.  By supposition $ \{x\} $ and $ \{y\} $ are $ \wedge_\tau $- sets i.e. $ \{x\}=\{x\}_\tau^\wedge$ and $ y \not \in \{x\}_\tau^\wedge $. Therefore there exists an open set $ V' $ such that $x \in V' $,  but $ y \not \in V' $. Similarly, since $ \{y\}=\{y\}_\tau^\wedge$,  there exists an open set $ U' $ such that $ y\in U' $ and  $ x \not \in U' $.  Hence $x , y$ are weakly separated by open sets $ V' $ and $ U'$  respectively
and $ (X,  \tau) $  is $ T_1 $-space.  
\end{proof}

\textbf{Definition  6.16} : A space $ (X,  \tau) $ is said to be  weak $ R_0$-space  if every $ \lambda^*$-closed   singleton is a $ \wedge_\tau $-set. \\

\textbf{Theorem 6.17 } : Every $ R_0 $-space  is a weak $ R_0 $-space.
\begin{proof}

Suppose $ x\in X $ and $ \{x\} $ is $ \lambda^* $-closed, then $ \{x\}=\{x\}_\tau^\wedge\cap\overline{\{x\}}$ by  Lemma 5.3 (iv). We claim that $\{x\} $ is a $ \wedge_\tau $-set. If not, then $  \{x\}\not=\{x\}_\tau^\wedge  $ and so let $  y\in \{x\}_\tau^\wedge -\{x\}$. Then $ y\not\in\overline{\{x\}}$.  So there is a closed set $F,   F\supset\{x\}  $  such that $ y\not\in F $. This implies that $ y\in X - F $,  an open set.  Since$ (X,\tau) $ is $ R_0  $-space,  $ \overline{\{y\}}\subset X - F$. Therefore  $\overline{\{y\}}\cap F=\emptyset$.  Since $ x\in F  $,  $ x\not\in\overline{\{y\}}=\{y\}\cup \{y\}'$ where $ \{y\}' $ denotes the set of limit points of $ \{y\} $. Therefore there exists an open set $ V \supset\{x\} $ such that $y\not\in V$,  since $x\not=y$ and $ x$ is not also the limit point of $ \{y\} $. This implies that $ y\not\in \{{x}\}_\tau^\wedge $,  a contradiction. Hence $ (X,  \tau) $ is weak $ R_0 $-space.
\end {proof}

But the following example shows that the converse of the Theorem 6.17 may  not be true.\\

\textbf{Example  6.18 :}  Let $ X=R $ and  $ G_i $  be the countable subsets of $X-Q-\{\sqrt{2}\}$, $ \tau$ = $\{X,  \emptyset,  G_i \}$. Then  $ (X,  \tau) $ is a space but not a topological space. We verify which singletons are $ \lambda^*$-closed.   Take $ \sqrt{2}\in X $.  Now $ \{\sqrt{2}\}_\tau^\wedge =X $ and $\{\overline{\sqrt{2}}\}=\{\sqrt{2}\}\cup Q $ which imply that $\{\sqrt{2}\} $ is not $ \lambda^*$-closed. Rational singletons are not $ \lambda^*$-closed. For, suppose that $A=\{\frac{1}{2}\}$. Then $ A_\tau^\wedge=X, \overline{A}=A \cup Q\cup \{\sqrt{2}\}=\{\sqrt{2}\}\cup Q $. Therefore $A_\tau^\wedge\cap \overline{A}=\{\sqrt{2}\}\cup Q \not=A$. Now for any irrationals $ r $ except $ \sqrt{2}, \{r\} $ is  $ \lambda^*$-closed. So $ (X, \tau) $ is a weak $ R_0 $- space.  Also if $ B = \{\sqrt{5}\} $ then $\overline{B}=B\cup Q\cup \{\sqrt{2}\}$ and  $B$ is an open set but $B$ does not contain $\overline{B}$, hence the space is not $ R_0 $.\\

\textbf{Lemma  6.19  :}  If every subset of $X$ is $ \wedge_\tau $-set, then $ (X, \tau) $ is $ T_1 $-space.

\begin{proof}
Let every subset of $X$ be $ \wedge_\tau $-set. Then every singleton is $ \wedge_\tau $-set. So by Theorem 6.15, $ (X, \tau) $ is $ T_1 $-space. 
\end {proof}

\textbf{Note:  6.20  :}  The converse is not true as revealed in the following Example 6.21. But the converse is true by imposing additional conditions as given in the Lemma 6.23. 
Also note that the converse part is true in a $\mu $-space  \cite{MS}.\\

\textbf{Example 6.21 :} Let $ X = R-Q $,    $ G_i $ be the countable subsets of $X$ and $ \tau $ = $\{X,  \emptyset, G_i\}$. Then $ (X, \tau) $ is a space but not a topological space. The space is $ T_1 $ by Theorem 6.15. Let $ B $ be the set of all irrationals in $ (0,  1) $. Then $ B $
is not a $ \wedge_\tau$-set as $ B_\tau^\wedge=X\not=B $.\\

\textbf{Definition  6.22}  (cf. Definition 15  \cite{DP}  : A space $ (X, \tau) $ is said to be strongly symmetric if  $ \{x\} $  is $ g^*$-closed  for each $x \in X$.\\

\textbf{Lemma  6.23  :}  If $ (X,  \tau) $ is a strongly symmetric $ T_1 $-space and satisfies the condition $ C=C^* $, then every subset of X is a $ \wedge_\tau $-set.
\begin{proof}

Let $ (X,  \tau) $ be a strongly symmetric $ T_1 $-space satisfying the condition $ C=C^* $ and let $ A\subset X$ , $ x\in X $ and $ x\not\in A $. Then, by definition, $ \{x\}$ is $ g^*$-closed  . Since $ (X,\tau) $ is $ T_1 $-space,  $ \{x\}  $ is a $ \wedge_\tau $-set by Theorem 6.15 and so a $ \lambda^*$-closed set.  Therefore $ \{x\} $ is a closed set by Theorem 5.8.  Therefore $ X - \{x\} $ is an open set containing A.  So $ A=\cap \{X - \{x\}, x\in X - A\} $ which implies $A$ is a  $ \wedge_\tau $-set.
\end{proof}

\textbf{Remark  6.24:} If the space  $ (X,  \tau) $ is a strongly symmetric $ T_1 $-space and satisfies the condition $ C=C^* $, then union and intersection of two $ \lambda^* $-closed sets are $ \lambda^* $-closed sets.\\

\textbf{Theorem 6.25}  {\cite{MS} : For a space $ (X, \tau)$ , the following statements are equivalent:

(1)   $ (X, \tau)$  is $ T_1 $ 

(2)   $ (X, \tau)$  is $ T_0 $ and $ R_0 $

 (3)  $ (X, \tau)$ is $ T_0 $  and weak $ R_0 $.
 
 \begin{proof}
 $(1)\Rightarrow (2)$: It follows from Theorem  6.14.
 
       $(2) \Rightarrow  (3)$: It follows from Theorem  6.17.
       
       For $(3)   \Rightarrow  (1)$:  let $ (X, \tau)$  be $ T_0 $  and weak $ R_0 $ and $ \{x\}\subset X $.  So by Theorem 6.4,  $ \{x\}$  is $ \lambda^*$-closed . Again  $ (X, \tau)$ is weak $ R_0 $,  $\{x\} $ is $ \wedge_\tau$-set. By Lemma 6.19,  $(X, \tau)$ is  $ T_1 $.
  \end{proof}
  
  \textbf{Theorem 6.26 :} If $ (X ,\tau) $ is a strongly symmetric $ T_1 $-space and satisfies the condition $ C = C^* $, then it is $ T\omega $-space.
  \begin{proof}
  Let $(X, \tau)$ be strongly symmetric $ T_1 $-space satisfying the condition $ C = C^*$ and let $ A $ be a $ g^*$-closed  set. Then $\overline{A^*}=A $ and so $\overline{A^*}$ is  $ g^*$-closed. Therefore $ A^c\in C^* = C $ which implies that $ \overline{A}$ is closed. Now let $ x\in\overline{A} - A $.  Since $ (X, \tau) $ is $ T_1  $ and strongly symmetric space, $ \{x\} $ is  $ g^*$-closed  and $ \lambda^*$-closed , and since $ C = C^*, \{x\} $ is a closed set by Theorem 5.8.  But by Theorem 3.5,  $\{x\}\not\subset \overline{A} - A$. Therefore $ x\in A $ and so $\overline{A}=A$ which implies that $A$  is a closed set and hence the space  $ (X, \tau) $ is $ T_\omega $.
  \end{proof}
  
  \textbf{Theorem 6.27  :} If  the space $(X, \tau) $ is strongly symmetric,  weak $ R_0 $,  and satisfies the condition $ C = C^* $, then the following are equivalent:
  
   (1)  $(X, \tau) $ is $ T_0 $
   
   (2)  $(X, \tau) $ is $ T_1 $
   
   (3)  $(X, \tau) $ is $ T_\omega $
   
   (4)  $(X, \tau) $ is $ T_\frac{5\omega}{8} $
   
   (5)  $(X, \tau) $ is $ T_\frac{3\omega}{8} $
   
   (6)   $(X, \tau) $ is $ T_\frac{\omega}{4} $.
   
   \begin{proof}
   
   (1) $  \Rightarrow(2) $: It follows from Theorem 6.25. 
   
   (2) $ \Rightarrow (3) $: It follows from Theorem 6.26.
   
   (3)  $\Rightarrow $ (4):  It follows from Theorem 5.6 and 6.9.
   
   (4)  $\Rightarrow $ (5):  It follows from Theorem 6.9 and 6.7.
   
   (5) $ \Rightarrow $(6) :  It follows from Theorem 6.7 and 6.3.
   
   (6)$ \Rightarrow $ (1): It follows from Theorem 6.3 and 6.2.   
  \end{proof}

\end{document}